\documentclass[11pt]{article}
\usepackage{latexsym}
\usepackage{amsfonts}
\usepackage{amsmath}
\makeatletter

\def\@footnotetext#1{\insert\footins{%

\footnotesize 

 \interlinepenalty\interfootnotelinepenalty

 \splittopskip\footnotesep

 \splitmaxdepth \dp\strutbox \floatingpenalty \@MM

 \hsize\columnwidth \@parboxrestore

 \edef\@currentlabel{\csname 
p@footnote\endcsname\@thefnmark}\@makefntext
 {\rule{\z@}{\footnotesep}\ignorespaces
#1\strut}}}

\def\abstract{\small\quotation{\hskip-\parindent\sc Abstract.}}
\def\classification{\@ifnextchar [{\@xfootnotenext}%
 {\begingroup\let\protect\noexpand
 \xdef\@thefnmark{}\endgroup
 \@footnotetext}}

\title {}

\begin{document}

\newcommand{\Spec}{\operatorname{Spec}}
\newcommand{\Aff}{{\mbox{\bf A}}}
\newcommand{\Mat}{\operatorname{Mat}}
\newcommand{\de}{\partial} 
\newcommand{\J}{\mbox{J}} 

\classification {{\it 2000 Mathematics Subject Classification:} Primary 
14E09, 14E25; Secondary 14A10, 13B25.
$\ast$)Partially supported by an NSA grant.\\
$\dagger$) Partially supported by Hong Kong RGC-CERG Grants
 10203186 and 10203669.}

\begin{center}

{\bf \Large  The Stable Equivalence and  Cancellation Problems}

\bigskip

{\bf Leonid Makar-Limanov$^{\ast}$ \hskip 20pt   Peter van Rossum \\

Vladimir Shpilrain \hskip 35pt  Jie-Tai Yu$^{\dagger}$}

\end{center}

\medskip

\begin{abstract}
\noindent    Let $K$ be an arbitrary field of characteristic 0, and 
 $\Aff^n$ the  $n$-dimensional affine space over $K$. 
A well-known cancellation problem asks,  given two 
algebraic varieties $V_1, V_2 \subseteq \Aff^n$ with isomorphic  
  cylinders $V_1 \times \Aff^1$ and  $V_2 \times \Aff^1$, whether  
  $V_1$ and  $V_2$ themselves are  isomorphic. 

In this paper, we  focus on  a related  problem: given two 
varieties with  {\it equivalent} (under an automorphism of $\Aff^{n+1}$) 
  cylinders $V_1 \times \Aff^1$ and  $V_2 \times \Aff^1$, are 
  $V_1$ and  $V_2$ equivalent under an automorphism of  $\Aff^n$? 
We call this {\it stable equivalence problem}. 
We show that the answer  is positive 
for any two curves  $V_1, V_2 \subseteq \Aff^2$. 

 For an arbitrary $n \ge 2$, we consider a special, arguably the 
most important,  case of both problems, where one of the varieties is
a hyperplane.  We show that a positive 
solution of the 
stable equivalence problem in this case implies a positive solution of 
the  cancellation problem. 

\end{abstract} 

\bigskip

\noindent {\bf 1. Introduction }

\bigskip

 Let $K[x_1,\dots, x_n]$ be the polynomial algebra in $n$ variables
over a field $K$ of characteristic 0. Any collection of polynomials 
$p_1,\dots,p_m$ from
$K[x_1,\dots, x_n]$ determines an algebraic variety
$\Spec K[x_1,\dots,x_n]/\langle p_1,\dots,p_m\rangle$ in the affine
space $\Aff^n = \Aff^n_K = \Spec K[x_1,\dots,x_n]$. 
If $K$ is algebraically closed and $\langle p_1, \dots, p_m\rangle$
is radical, we can of course think of this variety as the zero set
${\{}p_i=0, ~i=1,\dots,m{\}}$ in $K^n$.
We  denote this algebraic variety by $V(p_1,\dots,p_m)$.

We say that two algebraic varieties $V(p_1,\dots,p_m)$
 and $V(q_1,\dots,q_k)$ in $\Aff^n$  are
{\it equivalent} if there is an automorphism of $\Aff^n$ that takes one
of them onto the other. Algebraically, this means there is 
an automorphism of $K[x_1,\dots, x_n]$ that takes the ideal 
$\langle p_1,\dots,p_m \rangle$ to the ideal $\langle q_1,\dots,q_k 
\rangle$. 

 A variety equivalent to  $V \times \Aff^1$ is called a 
{\it cylinder};
 a variety of the form $V(p)$ is called a  {\it hypersurface}, and a 
hypersurface equivalent to $V(x_1)$ is called a  {\it hyperplane}.

We say that two algebraic varieties $V_1$ and  $V_2$ 
in $\Aff^n$ are {\it stably equivalent} 
if there is an automorphism of $\Aff^m$ for some $m > n$  
that takes  the cylinder $V_1 \times \Aff^{m-n}$
 onto  $V_2 \times \Aff^{m-n}$. 

We also say that  two 
polynomials $p, q \in K[x_1,\dots, x_n]$ are  
 stably  equivalent if $~\alpha(p) = q$ for some  
automorphism $~\alpha$ of $K[x_1,\dots, x_m]$, ~$m > n$.
 
 We address here the following 
\medskip

\noindent {\bf Stable equivalence problem.} Is it true that 
stable equivalence of 
two  hypersurfaces  in $\Aff^n$ implies their equivalence? 

 Or, in purely algebraic language: are any two stably equivalent 
polynomials equivalent?

\medskip

 If one considers arbitrary algebraic varieties, not just 
hypersurfaces, 
then the answer is negative, as explained in  \cite{ShYu2}. The 
corresponding 
example is based on a well-known example, due to Danielewski, of 
non-isomorphic
 surfaces in ${\bf C}^3$ with isomorphic cylinders.

 Here we solve the Stable equivalence problem for $n=2$: 
\medskip

\noindent {\bf Theorem 1.1.}  If two polynomials 
$p, q \in K[x, y]$ are stably equivalent, then they are equivalent. 
Or, in geometric language: if $V(p)$ and  $V(q)$ are two curves  in 
$\Aff^2$ 
such that, for some $s \ge 1$,  the cylinders $V(p) \times \Aff^{s}$ 
and 
$V(q) \times \Aff^{s}$ are equivalent in $\Aff^{2+s}$, then  $V(p)$ and  
$V(q)$ are 
equivalent in $\Aff^2$. 
\medskip 

  Upon replacing equivalence with isomorphism in the statement of 
Theorem 1.1, 
one gets a well known result of Abhyankar,  Eakin and  Heinzer 
\cite{AEH}.

 Now we focus on a special  case of the Stable equivalence problem; 
we call it 

\medskip

\noindent {\bf Stable coordinate conjecture.}  Let $V(p)$ be a 
hypersurface in $\Aff^n$. If $V(p) \times \Aff^1$ is equivalent to 
a hyperplane in $\Aff^{n+1}$, then $V(p)$ is equivalent to 
a hyperplane in $\Aff^n$. \\
 Or, in purely algebraic language: 
if $p=p(x_1,\dots, x_n)$ ~and   $~\varphi(p) = x_1~$ for some 
automorphism
 $~\varphi$ of $K[x_1,\dots, x_{n+1}]$, then  also $~\alpha(p) = x_1~$ 
for some automorphism $~\alpha$ of $K[x_1,\dots, x_n]$, i.e., 
$p$ is a {\it coordinate} in  $K[x_1,\dots, x_n]$. 
\medskip

 It turns out that the Stable coordinate conjecture is closely 
related to the famous Cancellation conjecture  of Zariski:

\medskip 

\noindent {\bf Cancellation conjecture}. 
   Let $V(p)$ be a hypersurface in $\Aff^n$.
If $V(p) \times \Aff^1$ is isomorphic to 
a hyperplane in $\Aff^{n+1}$ (i.e., to $\Aff^n$), then $V(p)$ is isomorphic to 
a hyperplane in $\Aff^n$ (i.e.,  to $\Aff^{n-1}$). \\
 Or, in purely algebraic language: if, for some $K$-algebra $R$,  
$R[x]$ is isomorphic to 
$K[x_1,\dots, x_{n}]$, then $R$ is isomorphic to $K[x_1,\dots, 
x_{n-1}]$. 
\medskip 

 This conjecture was proved for $n=2$ \cite{AEH}, \cite{Miy1} and  
$n=3$ \cite{MiySugie}, \cite{Fujita}. There is some circumstantial 
evidence that it 
might be wrong in  higher  dimensions if $K={\bf R}$, see  
\cite{Asanuma}. We refer to  \cite{Kraft} for a more detailed 
survey on this problem.

 In \cite{ShYu2}, it was shown that, for each 
particular $n$, the 
Cancellation conjecture follows from the Stable coordinate conjecture 
combined  with the Embedding conjecture of Abhyankar  and  Sathaye 
(see \cite{AM}), and also  that the Stable coordinate conjecture 
follows from the Cancellation conjecture combined  with the Embedding 
conjecture.

 Here we  establish a more straightforward implication: 

\medskip

\noindent {\bf Theorem  1.2.} For each 
particular $n$, the Stable coordinate conjecture  implies 
the Cancellation conjecture. 
\medskip

 It would be interesting to pinpoint also some connection between 
more general forms of both conjectures, namely, between what 
we call the Stable equivalence problem and  the Cancellation problem 
(see the abstract). In particular, having in mind Danielewski's example 
mentioned before and motivated by Theorem   1.2, we ask: 

\medskip

\noindent {\bf Problem 1.} Let $p=p(x,y,z)=xy+z^2$. Is it true that 
every  polynomial in $K[x,y,z]$  which is stably equivalent  to $p$ 
is, in fact,  equivalent  to $p$ ? 
\medskip

 Recall that, by  results of Danielewski \cite{Dan} and  Fieseler 
\cite{Fieseler}, the hypersurface 
$D(k)=\{xy^k+z^2+1=0\}$ is not isomorphic to $D(m)=\{xy^m+z^2+1=0\}$ 
 if $k \ne m, ~k, m \ge 1$,  
whereas the cylinders $D(k) \times {\bf C}$ and  
$D(m) \times {\bf C}$ are isomorphic.

 Finally, we mention that it would be also interesting to find any 
relation between the general Cancellation problem and  the general 
Embedding problem. A somewhat bold conjecture would be 
 that if, for a  hypersurface $V(p)\subseteq \Aff^n$, the cylinder 
 $V(p)\times \Aff^1$ has a unique (up to an automorphism 
of  $\Aff^{n+1}$) embedding into  $\Aff^{n+1}$, then, whenever $V(p) 
\times \Aff^1$ is isomorphic to $V(q) \times \Aff^1$, 
one has $V(p)$ isomorphic to $V(q)$. 

 Now a natural question is whether or not Danielewski's surfaces/cylinders  
have unique embeddings in ${\bf C}^4$. We were able to prove 
that all but one of them do not:
\medskip

\noindent {\bf Proposition 1.3.} For any $m \ge 2$, the hypersurface 
$D(m)\times {\bf C}^{k-3}=\{xy^m+z^2+1=0\}$ has at least 2  
inequivalent embeddings in ${\bf C}^k$ for any $k \ge 3$. 
\medskip

 We note that for $k = 3$, this was also proved in \cite{FM} 
(by an altogether different method). 
 For $m=1$, the question is open: 

\medskip

\noindent {\bf Problem 2.} Does the hypersurface  $D(1)\times {\bf 
C}=\{xy+z^2+1=0\}$ 
have a unique embedding in ${\bf C}^4$ ? 

\medskip 

 It is also unknown whether $D(1)$ has a unique embedding in ${\bf C}^3$ 
(cf. \cite[Question 1]{FM}). 

 We note that Problems 1 and 2 cannot both have positive  answers. 
 Indeed, if the answer to Problem 2 was positive, then, since we know 
that $D(1)\times {\bf C}$ is isomorphic to $D(m)\times {\bf C}$ 
for any $m \ge 1$, 
we would have $D(1)\times {\bf C}$ equivalent to $D(m)\times {\bf C}$ 
in ${\bf C}^4$. 
Then, if the answer to Problem 1 was positive, this would imply 
that $D(1)$ is equivalent to $D(m)$ in ${\bf C}^3$, which is known 
not to be the case. 

 This simple trick also works in a more general situation, namely: 
\medskip

\noindent {\bf Proposition 1.4.} Let $V(p), ~p=p(x_1,...,x_n)$
 be a hypersurface in $\Aff^n$. Suppose that the following two 
conditions hold: 

\smallskip 

\noindent {\bf (i)} $V(p) \times \Aff^1$ has a unique embedding in $\Aff^{n+1}$. 

\smallskip 

\noindent {\bf (ii)} If $V(p)$ is stably equivalent  to $V(q)$ 
for some $q=q(x_1,...,x_n)$, then $V(p)$ is equivalent  to $V(q)$.

\smallskip 

 Then, whenever $V(p) \times \Aff^1$  is isomorphic to $V(q) \times \Aff^1$,
  it follows that 
$V(p)$  and  $V(q)$ are isomorphic subvarieties of $\Aff^{n}$. 
\medskip

 The proof is obvious; we omit the details. Equally obvious is the following 
\medskip

\noindent {\bf Proposition 1.5.} Let $V(p), ~p=p(x_1,...,x_n)$
 be a hypersurface in $\Aff^n$. Suppose that the following two 
conditions hold: 

\smallskip 

\noindent {\bf (i)} $V(p)$ has a unique embedding in $\Aff^{n}$. 

\smallskip 

\noindent {\bf (ii)} If $V(p) \times \Aff^s$  and   $V(q) \times \Aff^s$  
are isomorphic subvarieties of  $\Aff^{n+s}$ for some $s \ge 1$,  
then $V(p)$  and  $V(q)$ are isomorphic subvarieties of $\Aff^{n}$. 
\smallskip 

 Then, whenever $V(p)$ is stably equivalent  to $V(q)$ 
for some $q=q(x_1,...,x_n)$, one has $V(p)$ equivalent  to $V(q)$.\\

\noindent {\bf 2. The two-variable case }

\bigskip

 Let $p, q \in K[x, y]$, ~$\varphi(p)=q$ 
for some automorphism $\varphi$ of $K[x, y, z, \dots]$. 

 Let $\varphi(x)=u=u(x,y,z,...), ~\varphi(y)=v=v(x,y,z,...)$. We are going to 
prove a stronger statement (Proposition 2.1 below) that will imply 
Theorem 1.1. 

We call a pair $(u,v)$ of polynomials {\it $z$-reduced} if the sum of  
$z$-degrees of the two polynomials cannot be reduced by either 
a (non-degenerate) linear transformation or a transformation of one of 
 the following two types: 

\noindent {\bf (i)} $(u, v) \longrightarrow (u+ \mu \cdot 
v^k, v)$ for some $\mu \in K^\ast; ~k \ge 2$;

\noindent {\bf (ii)}  $(u, v) \longrightarrow 
(u, v+ \mu \cdot u^k)$. 

 When proving Theorem 1.1, we can assume, without loss of generality, 
that the pair $(u(x,y,z,...), ~v(x,y,z,...))$ is  $z$-reduced. 
 
 If $(u,v)$ is a pair of {\it two-variable} polynomials such that 
the sum of their  degrees  
cannot be reduced by a transformation of one the above types, then we call 
this pair {\it elementary reduced}.

\medskip

\noindent {\bf Proposition 2.1.}  Let $p \in K[x, y]$ be a two-variable 
polynomial. 
Let $(u(x,y,z,...), ~v(x,y,z,...))$ 
be a $z$-reduced pair of algebraically independent polynomials    such 
that both of them actually depend on $z$. 
Then, for any ~$N \in {\bf Z}_+$, there is a polynomial  $w=w(x,y)$ 
such that 
deg($p(u(x,y,w, 0,...,0), v(x,y,w, 0,...,0))$) $> N$.  

\medskip 

 In the proof of Proposition 2.1,  we shall write just $u(x,y,z)$ 
 and  $v(x,y,z)$ instead of $u(x,y,z, 0,...,0)$ and  $v(x,y,z, 0,...,0)$ 
to simplify the notation.  
 First we prove 

\medskip

\noindent {\bf Lemma 2.2.}   Let $u(x,y, z)$ and  
$v(x,y,z)$ be algebraically independent. 
 For any  $M \in {\bf Z}_+$ and  $m, n > M$, there is 
 $c \in K$ such that 
$u(x,y, x^m y^n +c)$ and  $v(x,y,x^m y^n +c)$ are algebraically 
independent. 
\medskip

\noindent {\bf  Proof.}  Recall that polynomials $f_1, \dots, 
f_m \in
K[x_1,\dots,x_n]$ are algebraically dependent over $K$ if and only if 
the  Jacobian matrix $D(f_1,\dots,f_m)$ has rank smaller than $m$.

Assume, by way of contradiction, 
 that for all $c \in K$ the polynomials \\
$u(x,~y,~x^m y^n + c), ~v(x,~y, ~x^m y^n+ c) \in K[x,y]$ are 
algebraically dependent. This means that for all 
$c \in K$, the matrix
$$
  D(u(x,y,x^m y^n + c),~v(x,y,x^m y^n + c)) 
  =
  D(u,v) |_{(x,y,x^m y^n +c)} \cdot 
  \left(\begin{array}{cc}
    1 & 0 \\
    0 & 1 \\
    mx^{m-1} y^n & n x^m y^{n-1}
  \end{array}\right)
$$
has rank at most 1.  Then for all $c \in K$ the $2 \times 3$ 
matrix $D(u,v)|_{(x,y,x^m 
y^n+c)}$ has rank at most one, which means that all its $2 
\times 2$
minors are 0.  Using the fact that for all $a, b \in K$ 
the map from
$K$ to $K$ definied by $c \mapsto a^m b^n + c$ is surjective, this 
implies that
for all $a, b, c \in K$ all $2 \times 2$ minors of $D(u,v)|_{(a,b,c)}$ 
are 0.  Since $K$ is an infinite field, this in turn implies 
that all $2
\times 2$ minors of $D(u,v)$ are 0. Thus, 
 the rank of $D(u,v)$ is at most one and therefore $u$ and $v$ are 
algebraically dependent over $K$, a contradiction.  $\Box$
\medskip

\noindent {\bf  Proof of Proposition 2.1.} Since both $u(x,y, z)$ and  
$v(x,y,z)$ actually depend on $z$ and the pair is $z$-reduced, 
we can find $m, n > 2N$ and   $c \in K$ such that 
$u(x,y, x^m y^n +c)$ and  $v(x,y,x^m y^n +c)$ are algebraically 
independent (by Lemma 2.2) and elementary reduced. 

 Now we use a result of Shestakov     and  Umirbaev \cite{Umirbaev} 
which implies, in particular, that, if two  polynomials $r(x,y)$ and  
$s(x,y)$ of degree $> 2N$ are 
algebraically independent  and elementary reduced, then every non-constant 
 polynomial in the algebra $K[r,s]$   has degree  at least  
  $N+2$. This completes the proof of Proposition 2.1. $\Box$ 
\medskip

 Now we can get to the 
\medskip

\noindent {\bf  Proof of Theorem 1.1.} Recall that 
$\varphi(x)=u=u(x,y,z,...), 
~\varphi(y)=v=v(x,y,z,...)$. Upon applying an automorphism of $K[x,y]$
to both $p(x,y)$ and $q(x,y)$ if necessary, we may assume that 
$u(x,y,z,...)$ and  $v(x,y,z,...)$ are $z$-reduced. Now we have 
several cases.

\medskip

\noindent {\bf Case 1.} Both $u(x,y,z,...)$ and  $v(x,y,z,...)$ 
actually depend on $z$. Then we can apply Proposition 2.1 to get 
a contradiction  in this case. 
\medskip

\noindent {\bf Case 2.}   Say, $v(x,y,z,...)$ actually depends on $z$, 
whereas $u(x,y,z,...)$ does not.  
 Let  $x^m y^n$ be the highest monomial in $p(x,y)$ with respect to 
``lexdeg'' ordering with $y > x$. This monomial  will contain the 
highest power of $z$ after we plug in $u$ for $x$  
and  $v$ for $y$. This highest power of $z$ then cannot cancel out in 
$p(u,v)$. Therefore, $p(u,v)$ will depend on $z$, contrary to the 
assumption $p(u,v)=q(x,y)$. 
\medskip

\noindent {\bf Case 3.}  Neither  $u(x,y,z,...)$ nor 
$v(x,y,z,...)$ depend on $z$. If there are other variables that 
either $u(x,y,z,...)$ or 
$v(x,y,z,...)$ depend on, then we find ourselves in Case 1 or 2 above. 
If not, then there is
nothing to prove because   the restriction of $\varphi$ to $K[x,y]$ 
must be an automorphism of $K[x,y]$.

This completes the proof of Theorem 1.1. $\Box$  
\medskip

\noindent {\bf  Remark.} The crucial technical tool 
 in our proof of Theorem 1.1 was Shestakov-Umirbaev's result from  
\cite{Umirbaev} 
that bounds (from below) the degree of polynomials 
in the subalgebra of $K[x,y]$ generated by two given polynomials. 
This is (philosophically) similar to ``small cancellation'' ideas
in combinatorial  group  theory (see e.g. \cite{LS}). We note however 
that in commutative algebra, these ideas cannot be simply carried on
to higher  dimensions as the following example shows. 

Let $\varphi : x \to u=x-yt^2z^2, ~y \to v=1+tz^2, 
~z \to r=z^2, ~t \to s=-xt+yt^2+yt^3z^2$.  
Let  $p=p(x,y,z,t)= xy+zt$. 
 Then $\varphi(p)= x$, i.e., $x \in K[u,v,r,s]$,  even though 
the degrees of  $u,v,r,s$ are at least 2.  Similar examples can
be constructed with arbitrarily high  degrees of  $u,v,r,s$. 

 This example therefore makes it appear likely that our 
proof of Theorem 1.1 might  be difficult to carry on
to higher  dimensions, but, of course, this does not mean that 
the result itself does not hold.  
\medskip

 Finally, since our proof of Theorem 1.1 heavily relies on 
Shestakov-Umirbaev's result which is not yet published, 
we offer an alternative proof, which is  more elementary, 
but probably has more limited use. 
\medskip
 
\noindent {\bf Alternative proof of Theorem 1.1.}  The statement will 
follow from  Proposition 2.3 below.

\medskip

\noindent {\bf Proposition 2.3.}   Let $R = K[x_1, \dots, x_n, z]$
be  a ring of polynomials in $(n + 1)$ variables and let $(u, ~v)$
 be a $z$-reduced pair of algebraically
independent polynomials from $R$ such that both of them actually depend on
$z$. Then, for any nonconstant two-variable polynomial $p$,   
the  polynomial $p(u, v)$  depends on $z$, too.
\medskip

\noindent {\bf Proof.} 
Suppose, by way of contradiction, that 
$q = p(u, v) \in K[x_1, \dots, x_n]$.  It is clear that
$\frac{\de q}{\de x} \neq 0$ for some $x = x_i$: otherwise $q$ would 
be a constant, and $u$ and $v$  would be algebraically dependent, 
contrary to the  assumption.

Then the derivation $\de(g) = J_{x, z}(q, g) = \frac{\de q}{\de x}
\frac{\de g}{\de z} - \frac{\de q}{\de z} \frac{\de g}{\de x} =
\frac{\de q}{\de x}\frac{\de g}{\de z}$ is a nonzero
locally nilpotent derivation on $R$ , i.e., for any element $g \in R$ 
there is  $n$ such that $\de^n(g) = 0$. Indeed, since $\deg_z q = 0$, 
we see that $\deg_z \de (g) < \deg_z g$. Therefore, if 
$m = \deg_z g$, then $\de^{m+1}(g)= 0$.

 Define now a derivation $\de_1$ on $K[u, v] \subset
R$ as follows: $\de_1(f) = J_{u, v}(p(u, v), f)$ for $f=f(u, v) \in K[u, v]$. 
We claim that this derivation is locally nilpotent, too. 
Let $g=g(x_1, \dots, x_n, z) = f(u, v) \in K[u, v]$. By  the 
usual chain rule, we have $\de(g) = J_{x, z}(q, g) = 
J_{x, z}(p(u, v), f(u, v)) = J_{x, z}(u, v) \cdot J_{u, v}(p, f) = 
J_{x, z}(u, v)\cdot \de_1(f)$. Thus, 
if we consider $\de_1(f)$ as an element of $R$ 
(as opposed to just $K[u, v]$), then $\deg_z
\de_1(f) = \deg_z \de(g) - \deg_z J_{x, z}(u, v) < \deg_z g =
\deg_z f(u, v)$. As above, this implies that $\de_1$ is a locally
nilpotent derivation on $K[u, v]$ since every application of
$\de_1$ decreases the degree relative to $z$.

Locally nilpotent derivations on a polynomial ring in two
variables are well understood. In particular, it is known that the
kernel of a nonzero locally nilpotent derivation is a polynomial
ring in one variable and its generator is also a generator of the
ambient two-variable ring (see \cite{Ren}). Since $\de_1$ is a
nonzero derivation and $p \in \ker \de_1$, a generator $s$ of the
kernel  does not depend on $z$ either. 
Thus,  $K[u, v] = K[s, w]$, where
$s$ does not depend on $z$, and therefore 
$\deg_z u = \deg_w u \cdot  \deg_z w$,
$\deg_z v = \deg_w v \cdot  \deg_z w$. It is known (see e.g. \cite{Cohn}) 
that  $K[u, v] = K[y, w]$ implies 
that either $\deg_w u$ divides $\deg_w v$ or $\deg_w v$ divides
$\deg_w u$ and that there is an elementary transformation (see
the beginning of this section) which reduces the $w$-degree of the
pair $(u, ~v)$. Therefore, the $z$-degree of the pair can be reduced, too, 
so that  $(u, ~v)$ is not a $z$-reduced pair contrary to our assumption.
 This completes the proof. $\Box$\\

\noindent {\bf 3. The Stable  coordinate and other conjectures}

\bigskip

\noindent {\bf  Proof of Theorem 1.2.} One of the equivalent 
formulations of the Cancellation conjecture  is (see \cite[p. 
54]{Ebook}): 
for every locally nilpotent derivation $\cal D$ of the algebra 
$K[x_1,\dots, x_n]$ with a slice $s$, the kernel ~$\mbox{Ker}~\cal D$ is 
isomorphic to 
$K[x_1,\dots, x_{n-1}]$. By Proposition~2.1 of Wright \cite{Wright}, 
the latter 
property is equivalent to $s$ being a coordinate, i.e., an automorphic 
image of $x_1$. 

 Thus, we start with an arbitrary locally nilpotent derivation $\cal D$ 
of the algebra 
$K[x_1,\dots, x_n]$ with a slice $s$, and we want to prove that 
$s$ is a coordinate in $K[x_1,\dots, x_n]$. 

Extend $\cal D$ to 
$K[x_1,\dots, x_{n+1}]$ by ${\cal D}(x_{n+1})=0$. Then $\mbox{Ker}~\cal 
D$ in 
$K[x_1,\dots, x_{n+1}]$ is $K[x_1,\dots, x_n]^{\cal D}[x_{n+1}]$, where 
$K[x_1,\dots, x_n]^{\cal D}$ denotes ~$\mbox{Ker}~\cal D$ in 
$K[x_1,\dots, x_n]$. 
Since $s$ is transcendental over $K[x_1,\dots, x_n]^{\cal D}$, 
we have $K[x_1,\dots, x_n]^{\cal D}[x_{n+1}]$ isomorphic to 
$K[x_1,\dots, x_n]^{\cal D}[s]$. The latter algebra is equal to 
$K[x_1,\dots, x_n]$ by the result of Wright \cite[Proposition 2.1]{Wright}. 

 Thus we get  $K[x_1,\dots, x_{n+1}]^{\cal D}$ isomorphic to 
$K[x_1,\dots, x_n]$, which implies that $s$ is a coordinate in 
$K[x_1,\dots, x_{n+1}]$. Since we are under the assumption that the  
Stable  coordinate conjecture holds for this particular $n$, 
we conclude that $s$ is a coordinate in $K[x_1,\dots, x_n]$, 
and therefore the Cancellation conjecture holds for the same $n$. 
$\Box$ 

\medskip

\noindent {\bf  Proof of Proposition 1.3.} We give a proof here for 
$m=2$, just to simplify the notation. As in \cite{ShYu1}, it  will be 
technically more 
convenient  to write   algebras of residue classes as 
``algebras with relations'', i.e., for example, instead of 
$K[x_1,\dots,x_n]/\langle p(x_1,\dots,x_n) \rangle$ we shall write 
$\langle x_1,\dots,x_n \mid p(x_1,\dots,x_n)=0\rangle$. 

 We  get the following chain of ``elementary'' isomorphisms: 
\smallskip

\noindent 
  $\langle x, y, z ~\mid ~xy^2+z^2+1=0  \rangle ~\cong$ (applying 
the automorphism  $\phi : x \to x, ~y \to y+1, ~z \to z)$ \\
$\langle x, y, z ~\mid ~x= -xy^2-2xy-z^2-1 \rangle ~\cong \\
\langle x, y, z, u ~\mid ~u=xy, ~x= -uy-2u-z^2-1 \rangle ~\cong \\
\langle x, y, z, u ~\mid ~u=-uy^2-2uy-z^2y-y,  ~x= -uy-2u-z^2-1\rangle 
~\cong \\
\langle y, z, u ~\mid ~u=-uy^2-2uy-z^2y-y\rangle ~\cong \\
\langle x, y, z ~\mid ~x=-xy^2-2xy-z^2y-y\rangle ~\cong$  (applying 
the automorphism  $\phi : x \to x, \\
y \to y-1, ~z \to z)$ \\
$\langle x, y, z ~\mid ~xy^2+z^2y-z^2+y-1=0 \rangle$. 

\smallskip

 Now let $p=p(x,y,z)=xy^2+z^2+1, ~q=q(x,y,z)=xy^2+z^2y-z^2+y-1$. We are 
going 
to show that the gradients $grad(p)$ and  $grad(q)$ have different 
numbers of zeros. This obviously implies that $p$ and  $q$ are 
inequivalent 
under any automorphism of $K[x,y,z]$ (in fact, this implies that $p$ 
and  $q$ are 
even stably inequivalent). 

Compute: 
$$grad(p) = (y^2, ~2xy, ~2z)$$ 
$$grad(q) = (y^2, ~2xy+z^2+1, ~2yz-2z).$$

 We see that $grad(p)$ has infinitely many zeros ($y=z=0, ~x$ 
arbitrary), 
whereas $grad(q)$ has no zeros. This completes the proof. $\Box$ \\

\noindent {\bf Acknowledgements}
\medskip 

  The first 
and the third  authors are grateful to the 
 Department of Mathematics of the 
University of Hong Kong for its warm hospitality during their visit
when part of  this work was done. \\

\noindent 
 Department of Mathematics, Wayne State University,
 Detroit, MI  48202
\smallskip

\noindent {\it e-mail address\/}: 
lml@math.wayne.edu \\

\noindent 
 Department of Mathematical Sciences, New Mexico State University, Las 
Cruces,  NM 88011
\smallskip

\noindent {\it e-mail address\/}: 
petervr@nmsu.edu \\

\noindent 
 Department of Mathematics, The City  College  of New York, New York, 
NY 10031 
\smallskip

\noindent {\it e-mail address\/}: 
shpil@groups.sci.ccny.cuny.edu  

\smallskip

\noindent {\it http://www.sci.ccny.cuny.edu/\~\/shpil/} \\

\noindent Department of Mathematics, The University of Hong Kong, 
Pokfulam Road, Hong Kong 

\smallskip

\noindent {\it e-mail address\/}: 
yujt@hkusua.hku.hk
\smallskip

\noindent {\it http://hkumath.hku.hk/\~\/jtyu}

\end{document}